\documentclass{article}
\usepackage{amsfonts,amssymb,latexsym,amsmath,mathdots }
\newcounter{num}[section]
\newcommand{\Num}{\refstepcounter{num}%
\textbf{
\arabic{num}}}

\newcommand{\Theorem}{\textbf{Theorem}}
\newcommand{\Proof}{{\textbf{Proof}}}
\newcommand{\Def}{\textbf{Definition}}
\newcommand{\Conj}{\textbf{Conjecture}}
\newcommand{\Lemma}{ \textbf{Lemma}}

\newcommand{\Remark}{\textbf{Remark}}
\newcommand{\Prop}{\textbf{Proposition}}

\newcommand{\gx}{{\mathfrak g}}
\newcommand{\hx}{{\mathfrak h}}
\newcommand{\tx}{{\mathfrak t}}

\newcommand{\ux}{{\mathfrak u}}
\newcommand{\utxm}{{\ux\tx(m,\Fq)}}
\newcommand{\utxn}{{\ux\tx(n,\Fq)}}

\newcommand{\sx}{{\mathfrak s}}

\newcommand{\Gc}{{G^\circ}}
\newcommand{\gxd}{{\gx^\diamond}}
\newcommand{\Gd}{{G^\diamond}}

\newcommand{\Hc}{{H^\circ}}
\newcommand{\hxd}{{\hx^\diamond}}
\newcommand{\Hd}{{H^\diamond}}

\newcommand{\uxd}{{\ux^\diamond}}
\newcommand{\Ud}{{U^\diamond}}

\newcommand{\Kc}{{\mathcal K}}
\newcommand{\Xc}{{\mathcal X}}
\newcommand{\Dc}{{\mathcal D}}
\newcommand{\Oc}{{\mathcal O}}
\newcommand{\Ec}{{\mathcal E}}

\newcommand{\al}{{\alpha}}
\newcommand{\De}{{\Delta}}
\newcommand{\Dp}{{\De^+}}
\newcommand{\la}{{\lambda}}
\newcommand{\eps}{{\varepsilon}}

\newcommand{\UTm}{{\mathrm{UT}}(m,\Fq)}

\newcommand{\UT}{{\mathrm{UT}}}

\newcommand{\UTn}{{\mathrm{UT}}(n,\Fq)}
\newcommand{\utn}{{\ux{\frak t}(n,\Fq)}}
\newcommand{\Irr}{{\mathrm{Irr}}}

\newcommand{\Ch}{{\mathfrak{Ch}}}

\newcommand{\diag}{{\mathrm{diag}}}

\newcommand{\Fq}{{\mathbb F}_q}
\newcommand{\Cb}{{\mathbb C}}

\newcommand{\row}{{\mathrm{row}}}
\newcommand{\col}{{\mathrm{col}}}
\newcommand{\Ad}{{\mathrm{Ad}}}

\newcommand{\Ind}{{\mathrm{Ind}}}

\newcommand{\Matt}{{\mathrm{Mat}}}

\renewcommand{\leq}{\leqslant}

\begin{document}
\Large

\title{New supercharacter theory  for Sylow subgroups in  orthogonal and symplectic groups}
\author{A.N.Panov
\thanks{The classification of superclasses (section 2)  is supported by the RFBR grant  16-01-00154-a;
the construction of supercharacter theory (section 3)  is supported by  grant RSF-DFG grant  16-41-1013}}
\date{}
 \maketitle

\Large
\section*{\sc{\S 1  ~~Introduction. Supercharacters for unipotent groups}}

The notion of a supercharacter theory was suggested by P.Diaconis and I.M.Isaaks in the paper \cite{DI}.
 A priori every group affords several supercharacter theories.
 One of examples of a supercharacter theroy is the theory of irreducible characters.  Since for some groups
 (such as the unitriangular group, the Sylow subgroups in symplectic  and orthogonal groups, the parabolic subgroups and  others)  the problem of classification of irreducible characters (representations) remains to be a very complicated "wild" problem, it appears reasonable to
  replace this problem by the problem of construction of a supercharacter theory, which provides  the best approximation of  theory of  irreducible characters.

 Let us formulate the definition of a supercharacter theory following the paper  \cite{DI}.
  Let $G$ be a finite group, $1\in G$ be a unit element.  Let $\Ch = \{\chi_1, \ldots, \chi_N$\} be a system of complex characters (representations) of the group  $G$. \\
 \Def\Num\label{superch}. The system of characters  $\Ch$  determines a supercharacter theory of $G$ if there exists a partition $\Kc = \{K_1,\ldots, K_N\}$ of the group  $G$ satisfying the following conditions: \\
S1) the characters of $\Ch$ are pairwise disjoint  (orthogonal);\\
S2) ~ each character  $\chi_i$ are constant on each subset $K_j$;\\
S3) ~ $\{1\} \in \Kc$.
\\
Under this definition each character of $\Ch$  is  referred to as  a  {\it  supercharacter}, each subset of $\Kc$ a {\it  superclass}. Observe that the number of supercharacters is equal to the number of superclasses. The square table  $\{\chi_i(K_j)\}$ is called a  \emph{supercharacter table.}

For each supercharacter  $\chi_i$, consider its support  $X_i$ (the subset of all irreducible constituents of  $\chi_i$).  Observe that the condition S3)  of Definition  \ref{superch} may be replaced by  following condition:\\
 S3')  The system of subsets  $\Xc=\{X_1,\ldots, X_N\}$ is a partition  of the system of irreducible characters  $\Irr(G)$. Moreover,
 here each character   $\chi_i$ differs from the character  $\sigma_i=\sum_{\psi\in X_i} \psi(1)\psi$ by a constant factor (see \cite{DI,P1,P2}).

For the unitriangular group  $\UTm$, the suitable supercharacter theory was constructed in the series  of papers of C.A.M. Andr\'{e} \cite{A1,A2,A3}.
This theory was generated for the algebra groups by  P.Diaconis and I.M.Isaaks ~\cite{DI}.
By definition, an algebra group is a group of the form  $G=1+J$, where  $J$ is an associative  finite dimensional nilpotent algebra. Superclasses in the algebra group  $G$  are the equivalence  classes for the equivalence relation: $g\sim g'$, where $g=1+x$ and $g'=1+x'$, if there exist  $a,b\in G$ such that $x'=axb$. The similar relation is defined for  $J^*$: by definition,  $\la\sim\la'$ if there exist  $a,b\in G$ such that $\la'=a\la b$ (here  $a\la b(x)=\la(bxa)$).

Fix a nontrivial character  $t\to\eps^t$ of the additive group of the field $\Fq$  with values in the group of invertible elements of the field $\Cb$.  Supercharacters $\chi_\la$ of a given  algebra group are  the induced characters  from  linear characters  $$\xi_\la(1+x)=\eps^{\la(x)}$$
 of right stabilizers of  $\la\in J^*$. The sets of characters  $\{\chi_\la\}$  and classes $\{K(g)\}$, where $\la$ and  $g$ run through the set of representatives of equivalence classes of  $J^*$ and  $G$ respectively, give rise to a supercharacter theory of the algebra group  $G$. Supercharacters  $\chi_\la$  afford  the analog of A.A.Krillov formula (see \cite{DI, P2}):
 $$\chi_\la (1+x) = \frac{|G\la|}{|G\la G|}\sum_{\mu\in G\la G}\eps^{\mu(x)}.$$

An unipotent group is not an algebra group in general. The outlined method is not valid for unipotent groups.  In this paper, we propose
   the new approach which  can be applied for a large class of unipotent groups, hypothetically.  The application of this approach for the Sylow subgroups of  orthogonal and symplectic  groups enables   to construct the supercharacter theory (see Theorem  \ref{mainth}), which is a bit better then  the one suggested in the papers   \cite{AN-1, AN-2, AFN, Andrews}.

Let us present the content of this approach.
Let  $U$ be an unipotent group that is a semidirect product  $U=U_1U_0$ with the normal subgroup  $U_1$.
Suppose that  $U_0$ is  an algebra group, i.e.  $U_0=1+\ux_0$, where $\ux_0$ is an associated finite dimensional nilpotent algebra.  The Lie algebra   $\ux$ of the group $U$  is a direct sum of two subalgebras  $\ux=\ux_0 \oplus \ux_1$, where $\ux_0$ is an associated  algebra,  and $\ux_1=\mathrm{Lie}(U_1)$ is an ideal in $\ux$.
 Since $U_0$ is an algebra group, for any  $a\in U_0$ and $x_0\in \ux_0$, the elements  $\ell_a(x_0)=ax_0$ and $r_a(x_0)=x_0a$ also belongs to  $\ux_0$.
The left and right actions of  $U_0$ on $\ux_0$ can be extended to the actions on  $\ux$ as follows
\begin{equation}\label{lrAction}
 \begin{array}{l} \ell_a(x)= \ell_a(x_0)+ \Ad_a(x_1),\\
r_a(x) = r_a(x_0) + x_1, \end{array}
  \end{equation}
where $ a\in U_0$ and  $x=x_0+x_1$, ~ $x_0\in \ux_0$, ~ $x_1\in U_1$.
 Observe that  the left and right actions of the subgroup $U_0$ commute, and
\begin{equation}\label{adad}
\ell_ar_a^{-1}(x)=\Ad_a(x),
\end{equation}
where  $\Ad_a$ is the adjoint operator for  $a\in U_0$.
\\
\Def\Num\label{equiux}. Let $x, x'\in \ux$. The element $x$ is equivalent  to  $x'$ if  there exists a chain of transformations of forms\\
1) ~ $x\to \ell_a(x)$, where  $a\in U_0$,\\
2) ~ $x\to \Ad_u(x)$,  where $u\in U$, \\
that maps  $x$ to $x'$. \\
Because of  (\ref{adad})  we  may substitute $r_a$ for  $\ell_a$ in 1).

  Fix $\Ad$-invariant bijective map  $f:U\to \ux$, ~$f(1)=0$. As a map  $f$ we can take the logarithm  $\ln$ (it requires  strong restrictions of characteristic of the field;  see below Definition \ref{Spr} for Sylow subgroups in orthogonal and symplectic groups).
 Introduce the equivalence relation on  $U$ as follows.\\
 \Def\Num\label{equigr}. Two elements  $u_1$ and $u_2$ of the group  $U$ are equivalent  if the elements  $f(u_1)$ and $f(u_2)$ from  $\ux$ are equivalent in the sense of Definition \ref{equiux}.

Consider the equivalence classes  $\{K(u)\}$;  hypothetically they are superclasses for some supercharacter theory.

Let us define  the left and right actions of  $U_0$ on the dual space $\ux^*$ by the formulas
$$\ell^*_a\la(x)=\la(r_a(x)),$$
$$r^*_a\la(x)=\la(\ell_a(x)).$$
The equivalence relation on  $\ux^*$ is defined similarly to  \ref{equiux} for $\ux$.\\
\Def\Num\label{equiuxstar}. Let $\la, \la'\in \ux^*$. The  element  $\la$ is equivalent to  $\la'$ if there exists a chain of transformations of forms\\
1) ~ $\la\to \ell^*_a(\la)$, where  $a\in U_0$,\\
2) ~ $\la\to \Ad^*_u(\la)$,  where $u\in U$, \\
that maps  $x$ to $x'$. \\
As above  $\ell_a^* (r_a^*)^{-1} \la = \Ad_a^*\la$; in definition, we  may substitute $r^*_a$ for  $\ell^*_a$.
Denote by  $\Oc(\la)$ the equivalence class of  $\la\in \ux^*$. 

In this paper, we present the classification of equivalence classes in $\ux$, ~$U$ and $\ux^*$ for the Sylow subgroups in  orthogonal and symplectic groups (see Theorems \ref{thuxcl},~\ref{thuxclstar}, \ref{cl}).
\\
\Conj\Num\label{mconj}. There  exists a system of characters of a given  finite unipotent group  $U$  of the form
\begin{equation}\label{sumsuper}
\chi_{\la}(u) =  c(\la) \sum_{\mu\in \Oc(\la)} \eps^{\mu(f(u))}, ~~\mbox{where}~~ c(\la)\in \Cb,~ c(\la)\ne 0,
\end{equation}
such that along with  the partition of the group  $U$  into the classes   $\{K(u)\}$, where $\la$ and $u$ run through the sets of representatives of equivalence classes in $\ux^*$ and  $U$ respectively, give rise to a supercharacter theory of the group  $U$.\\
\Remark\Num\label{mainrem} (see \cite{Andrews}).    Observe that the formula  (\ref{sumsuper}) defines the system of orthogonal functions on  $U$ (since the characters  $\{\eps^{\la(x) }\}$ of the abelian group  $\ux$ are pairwise orthogonal). Easy to see that the functions  (\ref{sumsuper}) are constant on the classes  $K(u)$. From $f(1)=0$ it follows  $K(1) =1$.
So, the functions (\ref{sumsuper}) always fulfil the conditions  S1, S2, S3. The main problem is to prove   existence of constants  $c(\la)$  such that the formula (\ref{sumsuper}) defines a character of some representation of the group  $U$.

\section*{\sc{\S 2~~ Superclasses of the  orthogonal and  symplectic  groups}}

The unitriangular group   $G=\UTm$ consists of all upper triangular matrices of order  $m$  with ones on the diagonal and entries from the finite field  $\Fq$. Assume that  the characteristic of field  $p>2$. The Lie algebra of the unitriangular group  $\gx=\utxm$ consists of upper triangular matrices  with zeros on the diagonal.

Consider the matrices  $$I_n=\left(\begin{array}{ccc}0&\dots&1\\ \vdots&\iddots&\vdots\\ 1&\dots&0\end{array}\right)
\quad \mbox{and }\quad J_{2n}=\left(\begin{array}{cc} 0&I_n\\ -I_n&0\end{array}\right).$$
 Let $m$ denote the dimension of standard representation of Lie algebras of types $B_n$,~ $C_n$, and $D_n$. That is  $m=2n+1$ for $B_n$, and  $m=2n$ for $C_n$ and $D_n$.
The matrix algebra  $\mathrm{Mat}(m,\Fq) $ affords the involutive antiautomorphism  $X\to X^\dag$, where
$X^\dag= I_{m}X^tI_{m}$ for $B_n$ and $D_n$, and  $X^\dag= J_{2n}X^tJ_{2n}$ for $C_n$.

The standard Sylow subgroup  $U$ in orthogonal and symplectic group consists of all $g\in G$ obeying $g^\dag=g^{-1}$.  Respectively, its Lie algebra  $\ux=\{x\in \gx:~ x^\dag=-x\}$.

We denote by $X^\tau$ the matrix transposed to  $X$ with respect to the second diagonal.
The Lie algebra  $\ux$   for $C_n$ and $D_n$  consists of matrices of the form
\begin{equation}\label{uxcd}
\ux=\left\{\left( \begin{array}{cc} X_0&X_1\\0&-X_0^\tau\end{array}\right)\right\},
\end{equation}
where $X_0\in \utxn$, ~ $X_1^\tau=X_1$ for $C_n$ and   $X_1^\tau=-X_1$ for $D_n$.
The Lie algebra  $\ux$ is a sum of two subalgebras  $\ux=\ux_0+\ux_1$, where
$$\ux_0=\left\{\left( \begin{array}{cc} X_0&0\\0&-X_0^\tau\end{array}\right)\right\},
\quad\quad \ux_1=\left\{\left( \begin{array}{cc} 0&X_1\\0&0\end{array}\right)\right\}.$$
The subalgebra  $\ux_1$ is an ideal in $\ux$, and $\ux_0$ is isomorphic to  $\utn$, and, therefore,  it has a natural structure of associative algebra.

The group $U$ is a semidirect product  $U=U_1U_0$, where
\begin{equation} \label{uone}U_1=\left\{\begin{pmatrix} E&B\\0&E\end{pmatrix}\right\}\quad \mbox{and}\quad U_0=\left\{\begin{pmatrix} A&0\\0& (A^\tau)^{-1}\end{pmatrix}\right\},\end{equation}
$B^\tau=B$ for $C_n$ and  $B^\tau=-B$ for $D_n$, and $A\in \UTn$. The subgroup  $U_0$ is isomorphic to $\UTn$, and, therefore, it is an algebra group.

In the case  $B_n$, the Lie algebra  $\ux$  consists of matrices of the form
\begin{equation}\label{uxb}
\left\{\left( \begin{array}{ccc}
X_0&X_1&X_2\\
0&0&-X_1^\tau \\
0&0&- X_0^\tau
\end{array}\right)\right\}
\end{equation}
where $X_0\in \utxn$, ~ $X_1$ is a  $n\times 1$  column, ~ $X_2$ is a  $n\times n$ matrix and  $X_2^\tau=-X_2$.
As above the Lie algebra  $\ux$ is a sum of two subalgebras  $\ux=\ux_0+\ux_1$, where
$$
\ux_0 = \left\{\left( \begin{array}{ccc}
X_0&0&0\\
0&0&0 \\
0&0&- X_0^\tau
\end{array}\right)\right\},\quad
\ux_1 = \left\{\left( \begin{array}{ccc}
0&X_1&X_2\\
0&0&-X_1^\tau \\
0&0&0
\end{array}\right)\right\}.$$
 The subalgebra $\ux_1$ is an ideal in  $\ux$, and  $\ux_0$ is isomorphic to  $\utn$, and,  therefore, it has a natural structure of associative algebra.

The group  $U$ decomposes  $U=U_1U_0$, where
\begin{equation}\label{mbb} U_1=\left\{\begin{pmatrix} E&v&-\frac{1}{2}vv^\tau+B\\0&1& -v^\tau\\ 0&0&E \end{pmatrix}\right\}
\quad \mbox{and}\quad U_0=\left\{\begin{pmatrix} A&0&0\\0&1&0\\0&0& (A^\tau)^{-1}\end{pmatrix}\right\},\end{equation}
$B^\tau=-B$, and $v$ is a $n$-column.
The subgroup $U_0$ is isomorphic to  $\UTn$, and it is an algebra group.

 Let us define  the left and right actions of the subgroup  $U_0$ on $\ux$ following the formula (\ref{lrAction}).
For $C_n$ and $D_n$, ~ $a=\mathrm{diag}(A, (A^\tau)^{-1})$ and $x\in \ux$, we have
$$\ell_a(x) =  \ell_a(x)= \ell_a(x_0)+ \Ad_a(x_1) =
\begin{pmatrix} AX_0&AX_1A^\tau\\0&-X_0A^\tau\end{pmatrix}$$
$$r_a(x) =  r_a(x_0) + x_1 =  \begin{pmatrix} X_0A&X_1\\0&-A^\tau X^\tau\end{pmatrix}.$$

Observe
$$ \ell_a(x) = \begin{pmatrix} A&0\\0&E\end{pmatrix}\begin{pmatrix} X_0&X_1\\0&-X^\tau\end{pmatrix}
\begin{pmatrix} E&0\\0&A^\tau\end{pmatrix} = a_1xa_1^\dag, ~~\mbox{where}~~ a_1=\mathrm{diag}(A,E);$$
$$ r_a(x) = \begin{pmatrix} E&0\\0&A^\tau\end{pmatrix}\begin{pmatrix} X_0&X_1\\0&-X^\tau\end{pmatrix}
\begin{pmatrix} A&0\\0&E\end{pmatrix} =  a_2xa_2^\dag, ~~\mbox{where}~~ a_2=\mathrm{diag}(E,A^\tau).$$

For $B_n$, ~ $a=\mathrm{diag}(A,1,(A^\tau)^{-1})$, and $x\in \ux$, we have
$$\ell_a(x) =  \ell_a(x_0)+ \Ad_a(x_1) = \begin{pmatrix} AX_0&AX_1&AX_2A^\tau\\
0&0&-X_1^\tau A^\tau\\
0&0&- X_0^\tau A^\tau\end{pmatrix},$$

$$r_a(x) = r_a(x) =  r_a(x_0) + x_1 = \begin{pmatrix} X_0A &X_1&X_2\\
0&0&-X_1^\tau \\
0&0&- A^\tau X_0^\tau \end{pmatrix}.$$

Observe
$$ \ell_a(x) =
 \begin{pmatrix} A&0&0\\0&1&0\\ 0&0&E\end{pmatrix}
\begin{pmatrix} X_0&X_1&X_2\\
0&0&-X_1^\tau \\
0&0&- X_0^\tau\end{pmatrix} \begin{pmatrix} E&0&0\\0&1&0\\ 0&0&A^\tau\end{pmatrix} = a_1xa_1^\dag,$$  where $a_1=\mathrm{diag}(A,1,E);$
$$r_a(x) = \begin{pmatrix} E&0&0\\0&1&0\\ 0&0&A^\tau\end{pmatrix}
\begin{pmatrix} X_0&X_1&X_2\\
0&0&-X_1^\tau \\
0&0&- X_0^\tau\end{pmatrix} \begin{pmatrix} A&0&0\\0&1&0\\ 0&0&E \end{pmatrix}= a_2xa_2^\dag,$$
 where $a_2=\mathrm{diag}(E,1,A^\tau).$

Denote by  $\Gc$ the subgroup in $G=\UT(n,\Fq)$ generated by the subgroup  $U$ and the  matrices   $\diag(A_1,A_2)$ in the case  $C_n$ and $D_n$ (respectively, $\diag(A_1,1,A_2)$  in the case $B_n$).\\
{\bf Remark}. The elements  $x,x'\in\ux$ are equivalent if and only if there exists  $g\in\Gc$ such that  $x'=gxg^\dag$.

In referred to above papers  \cite{AN-1, AN-2, AFN, Andrews}, the equivalence relation is a bit coarser: $x\sim x'$ if there exists  $g\in G$ such that $x'=gxg^\dag$.

Denote by  $\Hc$ the subgroup in $G=\UT(n,\Fq)$ generated by the subgroup  $U_1$ and the matrices   $\diag(A_1,E)$ in the case  $C_n$ and $D_n$ (respectively,  $\diag(A_1,1,E)$ in the case  $B_n$).\\
Let us describe the subgroups $\Gc$  and  $\Hc$.\\
\\
\Prop\Num. 1) In the case  $B_n$,
\begin{equation}\label{Bnzero} \Gc=\left\{\begin{pmatrix}A_0&A_1&A_2\\0&A_3&A_4\\0&0&A_5\end{pmatrix}\right\}\quad \mbox{and}\quad \Hc=\left\{\begin{pmatrix}A_0&A_1&A_2\\0&A_3&A'_4\\0&0&E\end{pmatrix}\right\},
\end{equation}
  where $A_0,A_5\in \UT(n-1,\Fq)$, ~$A_1,A_2,A_4$ are an arbitrary matrices of sizes  $(n-1)\times 3$, ~$(n-1)\times (n-1)$, ~$3\times (n-1)$ respectively,~~ $A_4'$ is an arbitrary  $3\times n$ matrix with zero last row,~~ $E$ is the unit $n\times n$ matrix,  $A_3$ is the  $3\times 3$ matrix of the form
  \begin{equation}\label{expc}
 \begin{pmatrix}1&c&-\frac{1}{2}c^2\\0&1&-c\\ 0&0&1\end{pmatrix},\quad c\in \Fq
.\end{equation}
  2) In the case  $C_n$, the subgroup  $\Gc$ coincides with  $G$ and
     $ \Hc=\left\{\begin{pmatrix}A_0&A_1\\0&E\end{pmatrix}\right\}. $\\
    3) In the case $D_n$,
\begin{equation} \label{CDnzero}
\Gc=\left\{\begin{pmatrix}A_0&A_1\\0&A_2\end{pmatrix}\right\} \quad\mbox{and}\quad  \Hc=\left\{\begin{pmatrix}A_0&A_1\\0&E\end{pmatrix}\right\},
\end{equation}
 where $A_0,A_2\in \UTn$, and  $A_1$ is an arbitrary   $n\times n$ matrix of the form
 \begin{equation}\label{wd}
 \begin{pmatrix} *&*&\cdots&*\\
 	\vdots&\vdots&\ddots&\vdots\\
 	c&*&\cdots&*\\
 	0&-c&\cdots&*\end{pmatrix},\quad c\in \Fq.\end{equation}
 	\Proof. We shall prove for the subgroup  $\Hc$ (for $\Gc$ similarly). \\
\emph{Case $C_n$}.   The subgroup  $\Hc$ is generated by the matrices of the form \\ $\begin{pmatrix} A&0\\0&E\end{pmatrix}$, where $A\in \UTn$ ~ and  ~$\begin{pmatrix} E&B\\0&E\end{pmatrix}$ with $B^\tau = B$.
Then the matrix
 $$ \begin{pmatrix} A&0\\0&E\end{pmatrix}\begin{pmatrix} E&B\\0&E\end{pmatrix} \begin{pmatrix} A^{-1}&0\\0&E\end{pmatrix} =  \begin{pmatrix} E&AB\\0&E\end{pmatrix}$$
 also belongs to $\Hc$.
 Easy to verify that the linear subspace,  spanned by the matrices  of the form $AB$, where $A\in\UTn$ and $B ^\tau=B$, coincides with $\Matt(n,\Fq)$. This proves statement 2).\\
 \emph{ Case $D_n$}.  It is treated  similarly. Easy to verify  that
  the linear subspace,  spanned by the matrices of the form  $AB$, where  $A\in\UTn$ and $B^\tau=-B$, coincides with the subspace of matrices of the form  (\ref{wd}). This proves statement  3).\\
 \emph{Case $B_n$}.
 For any two  $n$-columns  $v_1$ and $v_2$ we consider the matrix
\begin{equation}\label{Mvv}
M(v_1,v_2) = \begin{pmatrix} E&v_1&-\frac{1}{2}v_1v_2^\tau\\0&1& -v_2^\tau\\ 0&0&E \end{pmatrix}.
\end{equation}
The group  $\Hc$ is generated by the matrices   $\diag(A,1,E)$ with $A\in \UTn$,  the matrices
$M(v,v)$, where $v$ is a $n$-column, and
 $$F(B) = \begin{pmatrix}
1&0&B\\0&1&0\\0&0&1
\end{pmatrix},$$  where $B^\tau =-B$. Analogically the case  $D_n$ one can show that  the matrices of the form  $F(B)$, where
$B$  is a matrix  (\ref{wd}), belong to  $\Hc$.

The equality
$$ \begin{pmatrix}
A&0&0\\ 0&1&0\\0&0&E
\end{pmatrix}
\cdot M(v,v)\cdot
 \begin{pmatrix}
A^{-1}&0&0\\ 0&1&0\\0&0&E
\end{pmatrix}
= M(Av, v)$$
implies the subgroup $\Hc$ contains all matrices of the form
$M(v_1,v_2)$  for all columns  $$v_1=\begin{pmatrix}
\beta_1\\ \vdots \\ \beta_n
\end{pmatrix}, \quad v_2=\begin{pmatrix}
\beta_1'\\ \vdots \\ \beta_n'
\end{pmatrix}~~~ \mbox{such~~ that}~~~ \beta_n=\beta_n'\ne 0. $$

For any  $n$-columns $v_1,v_2,d_1,d_2$ we have an equality
\begin{equation}\label{Mprod}
M(v_1,v_2) M(d_1,d_2)=  M(v_1+d_1,v_2+d_2) F(B_0),
\end{equation}
where $B_0=\frac{1}{2}(-v_1d_2^\tau + d_1v_2^\tau)$.

The equality  (\ref{Mprod}) implies the subgroup  $\Hc$ contains all matrices of the form  $F(B)$, where $B=(b_{ij})$ is an arbitrary  $n\times n$ matrix with $b_{n1}=0$. Finally, applying (\ref{Mprod}) one can verify $\Hc$ contains all matrices of the form  $M(v_1,v_2)$, where $\beta_n=\beta_n'$. This follows statement 1). $\Box$

Let us describe the  equivalence classes in  $\ux$ and $\ux^*$  (see Definitions  \ref{equiux} and \ref{equiuxstar}).
Order the set of integers of the segment  $[-n,n]$ as follows
$$1\prec\ldots \prec n\prec 0 \prec -n\prec \ldots \prec -1.$$
Denote by
 $\Delta^+$ the set of following integer pairs from  $[-n,n]$:\\
for $B_n$
$$ \Dp=\{(i,j):~~ 1\leq i\leq n,~~ i\prec j\prec -i\},$$
for $C_n$
$$ \Dp=\{(i,j):~~ 1\leq i\leq n,~~ i\prec j \preccurlyeq -i,~~ j\ne 0\},$$
for $D_n$
$$ \Dp=\{(i,j):~~ 1\leq i\leq n,~~ i\prec j\prec -i,~~ j\ne 0\}.$$
We  refer to elements from  $\Dp$ as  \textit{positive roots}, and to  $\Dp$ as the \textit{set of positive roots}.

For any positive root  $\al=(i,j)\in\Dp$, we call  $i$  a \textit{row number } (denote $i=\row(\al)$)
and $j$ a \textit{column number} (denote $j=\col(\al)$).
\\
\Def\Num. The subset $\Dc\subset \Dp$ is called  \textit{basic} if there is no more than one root from $\Dc$
in each row and each column.

The other name is the set of rook placement type. \\
\Def\Num. We refer to a subset  $\Dc\subset \Dp$ as  \textit{quasibasic } if
1) there is no more than one root of $\Dc$ in any column;\\
2) there is no more than one root of $\Dc$ in any row except the cases of  $B_n$,~  $D_n$ and pairs of roots  $(i,n)$ and  $(i,-n)$.

So, any quasibasic subset in  $C_n$ is a basic subset.
For any positive  root  $\al=(i,j)$ denote  $\al'=(-j,-i)$ (according to definition  $\al'$ is not a positive root). For any matrix  $x=(x_\al)\in \ux$ the entices  $x_\al$ and $x_{\al'}$ differs by a sign  $x_\al=\epsilon(\al)x_{\al'}$.

The Lie algebra  $\gx$ has the standard basis  $\{E_{ij}:~~ 1\leq i < j\leq m\}$. The Lie algebra  $\ux$ also has the standard basis  $\{\Ec_\al =E_\al+\epsilon(\al) E_{\al'}\}$.
By a pair  $(\Dc,\phi)$, where $\Dc$ is a quasibasic subset of $\Dp$ and a map $\phi: D\to \Fq^*$,
we define the element
$$  x_{\Dc,\phi} = \sum_{\al\in \Dc} \phi(\al) \Ec_\al.$$
\Theorem\Num\label{thuxcl}. 1) Each  element  $x\in \ux$ is equivalent to some element  $ x_{\Dc,\phi}$.
2)  The pair  $(\Dc,\phi)$ is uniquely determined by  $x$.
\\
\Proof. 1) Applying transformations $x\to gxg^\dag $, ~ $g\in \Gc$, we are able to obtain more zeros in the matrix $x$, and finally get $ x_{\Dc,\phi}$. \\
2)  Let us show that $(\Dc,\phi)$ is uniquely determined by  $x$.
Suppose that the elements  $x_{\Dc,\phi}$ and $x_{\Dc',\phi'}$ are equivalent.

For any positive root $(i,j)$, we consider the submatrix
$\Matt_{ij}(x)$ of $x$ with systems of rows and columns  $\{k:~ i\preccurlyeq k \preccurlyeq j\}$. Easy to see that if $x\sim x'$, then  the submatrices
$\Matt_{ij}(x)$ and $\Matt_{ij}(x')$ have equal ranks.

Suppose that   $\Dc$ and $\Dc'$ are basic subsets. The equality of ranks implies that  $\Dc=\Dc'$.
  For each root  $\al=(i,j)\in \Dc$, we consider the subset $\Dc_\al\subset \Dc$ that consists of
 $(k,m)\in\Dc$, where $k\succcurlyeq i$ and $m \preccurlyeq j$.  For $\al\in\Dc$ we consider the minor $M_{\al}$  of the matrix $x$ with systems of rows
$\row(\Dc_\al)$ and columns  $\col(\Dc_\al)$. It is not difficult to show that the equivalence  of matrices  $x_{\Dc,\phi}$ and $x_{\Dc,\phi'}$ implies  $M_\al(x_{\Dc,\phi}) = M_\al(x_{\Dc,\phi'})$. Hence $\phi=\phi'$.

Assume that the quasibasic subset  $\Dc$ is not basic. In this case, $\ux$ is of the type  $B_n$ or $D_n$,  and $\Dc$ contains the pair of roots   $\beta_1=(i,n)$ and  $\beta_2=(i,-n)$.
For $D_n$  the statement can be proves similarly the case of basic subset.

Consider the case of  $B_n$. The equality of ranks of matrices  $\Matt_{ij}$ implies that
$\Dc\setminus \{\beta_2\} = \Dc'\setminus \{\beta_2\}$.
  For $\al\in \Dc\setminus\{ \beta_2\}$, we consider the system of roots  $\Dc_\al\in\Dc$ that consists of   $(k,m)\in\Dc\setminus \{ \beta_2\}$, where $k\succcurlyeq i$ and $m \preccurlyeq j$. By the subset $\Dc_\al$ we define the minor   $M_\al $ as above. Since $M_\al(x_{\Dc,\phi}) = M_\al(x_{\Dc,\phi'})$, it follows   $\phi(\al)=\phi'(\al)$ for each  $\al\in \Dc\setminus \{ \beta_2\}$.

It remains to show  $\beta_2\in \Dc'$ and  $\phi(\beta_2)=\phi'(\beta_2)$.
For $\beta_1$, we construct the root system  $\widetilde{\Dc}_{\beta_1}$ that coincides with $\Dc_{\beta_1}$  if there is no root $(k,0)$, ~ $i<k\leq n$, in $\Dc$; if such root exists,  then  $\widetilde{\Dc}_{\beta_1} = \Dc_{\beta_1}\cup {(k,0)}$. As above we define the minor   $\widetilde{M}_{\beta_1}$.

For   $\beta_0=(i,0)$, we construct the system of roots   $\widetilde{\Dc}_{\beta_0}$ that coincides with  $$(\Dc_{\beta_1}\setminus\{\beta_1\})\cup \{\beta_0\}$$
if there is no root  $(k,0)$, ~ $i<k\leq n$, in  $\Dc$; if such root exists, then $\widetilde{\Dc}_{\beta_0} = \Dc_{\beta_1}\cup \{(k,-n)\}$.  As above we define the minor $\widetilde{M}_{\beta_0}$.

  Take $\widetilde{\Dc}_{\beta_2} = \Dc_{\beta_2}\setminus\{\beta_1\}$ and consider the corresponding minor
$\widetilde{M}_{\beta_2}$.
One can show that polynomial
$$I= \widetilde{M}_{\beta_1} \widetilde{M}_{\beta_2} + \frac{1}{2}
\widetilde{M}^2_{\beta_0}$$
is constant on the equivalence class of element  $x=x_{\Dc,\phi}$. As  $x\sim x'$, where  $x'=x_{\Dc',\phi'}$, we have  $I(x)=I(x')$.
Observe that  $\widetilde{M}_{\beta_0}(x)=\widetilde{M}_{\beta_0}(x')=0$, and $\widetilde{M}_{\beta_1}(x) = \widetilde{M}_{\beta_1}(x') \ne 0$. Therefore,   $\widetilde{M}_{\beta_2}(x) = \widetilde{M}_{\beta_2}(x') \ne 0$.
Hence  $\beta_2\in \Dc'$, the values of  $\phi$ and  $\phi'$ on the root $\beta_2$ coincide. ~$\Box$

The dual space  $\gx^*$  has the dual basis  $\{E^*_{ij}:~~ 1\leq i < j\leq m\}$.
The dual space  $\ux^*$  has the basis  $\{E^*_\al+\epsilon(\al) E^*_{\al'}\}$.
For  $(\Dc,\phi)$, we construct  the element
$$  \la_{\Dc,\phi} = \sum_{\al\in \Dc} \phi(\al) \left(E^*_\al+\epsilon(\al) E^*_{\al'}\right).$$
\Theorem\Num\label{thuxclstar}. 1) Each element  $\la\in \ux$ is equivalent to some element  $ \la_{\Dc,\phi}$.
2)  The pair $(\Dc,\phi)$ is uniquely determined by $\la$.\\
\Proof. The proof is similar to the one of Theorem \ref{thuxcl}.~ $\Box$

Turn to the equivalence relation on the  group $U$.
We take the Springer map as  a map  $f$ in Definition  \ref{equigr}.  \\
\Def\Num\label{Spr}. A map  $f:G\to\gx$ is called Springer map if $f$ is a bijection and  it obeys the following conditions:\\
1)  ~ $f(U)=\ux$,\\
2) ~ there exist $a_2,a_3,\ldots $ from the field $\Fq$ such that  $f(1+x) = x+a_2x^2+a_3x^3+\ldots$ for any  $x\in\gx$.

Examples of  Springer map are as follows:\\
1) the logarithm map  $\ln(1+x) = \sum_{i=1}^\infty (-1)^{i+1}\frac{x^i}{i}$~~ (it requires strong restrictions on the  characteristic of  field);\\
2) Cayley's map $f(1+x)=\frac{2x}{x+2}$ ~~(for $\mathrm{char}\,\Fq\ne 2)$.

 Denote   $u_{\Dc,\phi}=f^{-1}(x_{\Dc,\phi})$. Classification of equivalence classes in group $U$ follows from  Theorem \ref{thuxcl}.\\
\Theorem\Num\label{cl}. 1) Each element  $u\in U$ is equivalent to some  $ u_{\Dc,\phi}$.
2)  The pair $(\Dc,\phi)$ is uniquely determined by  $u$.

\section*{\sc{\S 3~~ Supercharacters of the  orthogonal and symplectic groups}}

In the case  $C_n$ and $D_n$, the group  $\Gc$  is  an algebra group, and construction of supercharacters for the group $U$ doesn't differ from the approach of paper  \cite{Andrews}.

In the case $B_n$, the group  $\Gc$ is not an algebra group.
Consider the subgroup  $S$ of all matrices  $\mathrm{diag}(E,M,E)$, where $M$ is a matrix of the type  (\ref{expc}).
The group  $\Gc$ is a semidirect product $\Gc=S\Gd$, where
$\Gd$ is a subgroup of $\Gc$ that consists of matrices of the form  (\ref{Bnzero}) with $A_3=E$.
The subgroup  $\Gd$ is an algebra group  $\Gd=1+\gxd$, where $\gxd$ is the Lie subalgebra of the group $\Gd$.

The group $U$ is also a semidirect product  $U=S\Ud$, where  $\Ud=U\cap \Gd$.
The Lie algebra  $\ux$ decomposes  $\ux=\sx\oplus\uxd$, where $\sx$ is the one dimensional Lie algebra spanned by the vector $\Ec_{n,0}$.

Take  $\Hd=\Gd\cap \Hc$. The subgroup  $\Hd$ is an algebra group  and  $\Hc=S\Hd$.
Easy to see that in each case  $B_n,~C_n,~D_n$ we have  $\Gc=U\Hd$ (see \cite[Lemma 6.6]{Andrews}).

For any linear form  $\eta\in (\gxd)^*$, we define the following associative subalgebras in  $\gxd$:\\
1) ~$r_\eta^\diamond=\{x\in\gxd:~ \eta(xy)=0~ ~\mbox{for~~any}~~ y\in \hxd\}$;\\
2) ~$\ell_\eta^\diamond=\{x\in\gxd:~ \eta(y^\dag x)=0~ ~\mbox{for~~any}~~ y\in \hxd\}$;\\
3) ~$\gx_\eta^\diamond= r_\eta^\diamond\cap \ell_\eta^\diamond$.

 Since $(r_\eta^\diamond)\dag = \ell_\eta^\diamond$, we have  $(\gx_\eta^\diamond)^\dag = \gx_\eta^\diamond$.  The subgroup  $\Gc$ contains the algebra subgroups  $R_\eta^\diamond=1+r_\eta^\diamond$, ~ $L_\eta^\diamond=1+\ell_\eta^\diamond$ and
 $G_\eta^\diamond=1+\gx_\eta^\diamond$.\\
 \Lemma\Num\label{and}~\cite[Lemma 6.3]{Andrews}. $\eta(xy)=0$ for any  $x,~y\in \gx_\eta^\diamond$. \\
 \Proof. Present  $x$ in the form  $x=x'+x''$, where $x'_{ij}=0$ for $1\preccurlyeq j\preccurlyeq 0$, and
 $x''_{ij}=0$ for $0\prec j\preccurlyeq -1$.
 Present $y$  in the form $y=y'+y''$, where  $y'_{ij}=0$ for $0\prec i\preccurlyeq -1$, and
 $y''_{ij}=0$ for $1\preccurlyeq i \preccurlyeq 0$. Then $x'\in(\hxd)^\dag$,~ $y'\in \hxd$, ~$x'y'=0$, ~ $x''y''=0$. The equality  $xy=x'y+xy'$ implies  $\eta (xy)=0$.  ~$\Box$

Let
 $\la\in(\uxd)^*$ and $\eta\in\gxd$ be such that  $\eta^\dag=-\eta$ and the restriction of $\eta$ on  $\uxd$ coincides with $\la$. Define   $$\ux_\la^\diamond = \ux\cap \gx_\eta^\diamond ~~\mbox{and}~~
U_\la^\diamond = U\cap G_\eta^\diamond.$$

For  $x\in \uxd $, we have
\begin{equation}\label{yxyx}
\eta(y^\dag x) =
-\eta^\dag (y^\dag x) = -\eta ((y^\dag x)^\dag) = -\eta(x^\dag y)=\eta(xy).
\end{equation}
 It follows    $ \ux\cap r_\eta^\diamond =  \ux\cap \ell_\eta^\diamond = \ux_\la^\diamond $.

 By Lemma  \ref{and}, the restriction of $\la$ on  the Lie algebra $\ux_\la^\diamond$ is its character with values in the field  $\Fq$.

 Let $\pi^\diamond$ be the natural projection  $(\uxd)^*\to (\ux_\la^\diamond)^*$. The following statement follows from the paper  \cite{Andrews}. For readers convenience we present it with  complete proof. \\
\Lemma\Num\label{pid}. For any  $\la\in(\uxd)^*$ the fiber  $(\pi^\diamond)^{-1} \pi^\diamond(\la)$  equals to $ \Hd\centerdot\la$.\\
\Proof.
\emph{Item 1}. Let $\eta\in (\gxd)^*$ and  $P$ is the natural projection $(\gxd)^*\to (r_\eta^\diamond)^*$.
Let us show that  $P^{-1}P(\eta)=\Hd\eta$ (here $h\eta(x)=\eta(xh)$ is the left action   $h\in\Hd$ on $\eta$).
Indeed, $P^{-1}P(\eta)=\eta + (r_\eta^\diamond)^\perp$.
 The definition of  $r_\eta^\diamond$ implies  $$r_\eta^\diamond=\{x\in \gxd:~~ y\eta(x)=0 ~ \mbox{for~any}~ y\in \hxd\} = (\hxd\eta)^\perp.$$
Hence $(r_\eta^\diamond)^\perp =\hxd\eta$ and $P^{-1}P(\eta)=\eta + \hxd\eta = (1+\hxd)\eta = \Hd\eta$.\\
\emph{Item 2}. Let $\Pi$ be the natural projection  $(\gxd)^*\to (\uxd)^*$. Let $\eta^\dag=-\eta$ and $\la= \Pi(\eta)$. Let us show that  $\Pi(\Hd\eta)=\Hd\centerdot \la$.

Really, for $h=1+y\in \Hd$ and $x\in\uxd$, we obtain
$h\centerdot \la(x) = (h\la h^\dag)(x) = \la(h^\dag xh) =
\la((1+y)^\dag x(1+y)) = \la(x) + \la(y^\dag x+xy) + \la(y^\dag x y). $

Observe that   $y^\dag x y=0$ for any  $x\in \gx$ and $y\in \hxd$.  Applying the equality (\ref{yxyx}) we obtain
$h\centerdot \la(x) = \eta(x)+2\eta(xy) = (1+2y)\eta(x).$ Hence $\Hd\centerdot \la = \Pi(\Hd\eta)$; this proves statement  2.\\
\emph{Item 3}. Since  $\ux_\la^\diamond= \uxd\cap r_\la^\diamond$,  we have
$$ (\pi^\diamond)^{-1} \pi^\diamond(\la) = \Pi(P^{-1}P(\eta)) = \Pi(\Hd\eta) = \Hd\centerdot \la. ~~\Box $$
  Let $\la=\la_{\Dc,\phi}$ and  $\eta=\eta_{\Dc,\phi}$ is the element $(\gx)^*$ such that $\Pi(\eta)=\la$ and $\eta^\dag=-\eta$.\\
 \Lemma\Num\label{lemlem}. Let  $\la$ and  $\eta$ are defined by  $\Dc,\phi$ as above.
 Then $\sx\gx_\la^\diamond\subset \gx_\la^\diamond$ and  $\gx_\la^\diamond\sx\subset \gx_\la^\diamond$.\\
 \Proof. Let $\eta\in (\gxd)^*$ and $\Pi(\eta)=\la$. For each $\al=(i,j)\in \Dp$, where $1\leq i<n$ and $i\prec j$,  we define the subalgebra  $r_\al^\diamond$ that consists of matrices  $x\in \gxd$ obeying the conditions:\\
 1) if $1\leq j<n$, then $x_{ik}=0$ for all  $i<k<j$;\\
 2) if $j\in\{n,0,-n\}$, then $x_{ik}=0$ for all   $i<k<n$;\\
 3) if $ -n<j\leq -1$, then $x_{ik}=0$ for all   $i\prec k \prec -n$, and $x_{-j,k}=0$ for all    $-j\prec k \prec -n$.

 Easy to see that the subalgebra  $r_\al^\diamond$ is invariant with respect to the left and right multiplication by  $\sx$.
The subalgebra   $r_\eta^\diamond$ coincides with intersection of subalgebras $r_\al^\diamond$ over all  $\al\in \Dc$.
 The  subalgebra   $r_\eta^\diamond$ is also invariant with respect to the left and right multiplication by  $\sx$.  $\Box$

Define the subalgebra  $\ux_\la$ in each of the following cases separately. \\
1) ~ $\Dc$ doesn't contain any roots of the form  $(i,0)$ and $(i,n)$.    In this case, we take  $\gx_\la=\sx\oplus\gx_\la^\diamond$ and $G_\la=SG_\la^\diamond$. Denote $\ux_\la = \gx_\la\cap \ux = \sx\oplus  \ux_\la^\diamond$. Take  $U_\la = G_\la\cap U =SU_\la^\diamond$.\\
2) ~$\Dc$  contains  $(i,0)$ or $(i,n)$.
In this case, we define  $\gx_\la=\gx_\la^\diamond$, ~ $\ux_\la=\ux_\la^\diamond$, and $G_\la=U_\la^\diamond$.

The associative subalgebra  $\hxd$ is an ideal in the associative algebra  $\gx$. Therefore,
$\hx_\la=\gx_\la + \hx^\diamond$ is its associative subalgebra  in $\gx$.
Then the subgroup  $H_\la=G_\la H = 1+\hx_\la$ is an algebra group in  $G$.

Let  $\pi$ be the natural projection  $\ux^*\to \ux_\la^*$.\\
\Lemma\Num\label{pipi}. 1) For any  $\la=\la_{\Dc,\phi}$,  the fiber $\pi^{-1}\pi(\la)$ coincides with  $ H_\la\centerdot\la$.\\
2) The formula
\begin{equation}\label{xixi}
\xi_\la(u)=\eps^{\la(f(u))}
\end{equation}
defines a character of the subgroup  $U_\la$.\\
\Proof.
\emph{Item 1.} ~$\Dc$ doesn't contain any roots of the form  $(i,0)$ and $(i,n)$.  Then $\eta(\sx x)=\eta(x\sx)=0$ for any  $x\in \gx$.  From Lemma  \ref{and} it follows that $\xi_\la(g)=\eps^{\la(f(g))}$ is a character of  $G_\la$; this proves statement 2).

By direct calculations, we obtain   $h\eta h^\dag(x)=\eta(x)$ for any
$x\in \ux_\la$ and  $h\in H_\la$. The Lemma \ref{pid} implies statement 1).\\
\emph{Item 2.}~ $\Dc$  contains  one of roots of type $(i,0)$ or $(i,n)$. Denote this root from $\Dc$ by $\gamma$.  If $\gamma = (i,0)$, then we take $z=E_{-n,-i} $; if $\gamma=(i,-n)$, then we take $z=E_{0,-i}$. In each case  $z\in \gx_\la\subset\hx_\la$ and   $\eta(\Ec_{n0}z)=\pm\phi(\gamma)\ne 0$.
The element  $h_t=1+tz$ belongs to $H_\la$ for any  $t\in \Fq$. By direct calculations, we obtain  $h_t\eta h^\dag_t(x)=\eta(x)$ for any  $x\in\ux^\diamond$,
      and  $h_t\eta h^\dag_t (\Ec_{n0}) = \eta(\Ec_{n0})\pm 2t\phi(\gamma)$. It follows  $\Pi^{-1}\Pi(\eta)$ is contained in $H_\la\centerdot \la$.
The Lemma  \ref{pid} implies statement 1). Statement  2) follows from Lemma \ref{and}.~
$\Box$

By the linear form  $\la=\la_{D,\phi}$, we define  a linear character of the subgroup  $U_\la$ as follows
$$\xi_\la(u)=\eps^{\la(f(u))}. $$
Consider the induced character
$$\chi_\la = \Ind(\xi_\la,U_\la, U).$$
\Theorem\Num\label{kirprop}. If  $\la=\la_{D,\phi}$, then $$\chi_\la(u) = \frac{|H_\la\centerdot\la|}{|\Gc\centerdot\la|}\sum_{\mu\in \Gc\centerdot\la} \eps^{\mu(f(u))}.$$
\Proof. Let  $\dot{\xi_\la}$ define the function on the group $U$ equal to  $\xi_\la$ on $U_\la$ and zero outside of  $U_\la$.
By definition
 $$
\chi_\la(u)=\frac{1}{|U_\la|}\sum_{v\in U} \dot{\xi_\la}(vuv^{-1}).$$
Applying Lemma  \ref{pipi} for all  $x\in \ux$,  we have
$$\sum_{\mu\in \pi^{-1}\pi(\la)} \eps^{\mu(x)} = \eps^{\la(x)}\cdot \sum_{\nu\in \ux_\la^\perp} \eps^{\nu(x)} =\left\{\begin{array}{l} \frac{|U|}{|U_\la|}\cdot\eps^{\la(x)}, ~~\mbox{if}~~ x\in \ux_\la;\\
 0, ~~\mbox{if}~~ x\notin \ux_\la.\end{array}\right.$$
 Hence  $$\dot{\xi}_\la(u)=\frac{|U_\la|}{|U|}\sum_{\mu\in \pi^{-1}\pi(\la)} \eps^{\mu(f(u))} = \frac{|U_\la|}{|U|}\sum_{\mu\in H_\la\centerdot\la} \eps^{\mu(f(u))}.$$
 Then
   \begin{multline*}
 \chi_\la(u)=\frac{1}{|U|}\cdot \sum_{\mu\in H_\la\centerdot\la,~v\in U} \eps^{\mu(vf(u)v^{-1})} = \frac{|H_\la\centerdot\la|}{|U|\cdot|H_\la|}\sum_{h\in H_\la,~ v\in U}\eps^{h\centerdot \la(vf(u)v^{-1})}  =\\
   \frac{|H_\la\centerdot\la|}{|U|\cdot|H_\la|}\sum_{h\in H_\la,~ v\in U}\eps^{\la(hvf(u)v^\dag h^\dag)}.
                                                                                                     \end{multline*}

Since  $\Gc=UH_\la$, we get
 \begin{multline*}
   \chi_\la(u) = \frac{|H_\la\centerdot\la|\cdot |U\cap H_\la|}{|U|\cdot|H_\la|}\sum_{g\in\Gc}\eps^{\la(gf(u)g^\dag)} = \frac{|H_\la\centerdot\la|}{|\Gc|}\sum_{g\in\Gc}\eps^{\la(gf(u)g^\dag)} = \\ \frac{|H_\la\centerdot\la|}{|\Gc\centerdot\la|}\sum_{\mu\in \Gc\centerdot\la} \eps^{\mu(f(u))}.  ~~\Box                                                                                        \end{multline*}
  \Theorem\Num\label{mainth}. Let  $U$ be the Sylow subgroup in  orthogonal or symplectic group.
The system of characters  $\{\chi_\la\}$,  and the partition of the group  $U$  into classes   $\{K(u)\}$, where $\la$ and  $u$ run through the sets of representatives of equivalence classes  $\la_{\Dc,\phi}\in \ux^*$ and  $u_{\Dc,\phi}\in U$, give rise to a supercharacter theory of the group  $U$.\\
\Proof. The proof follows from Remark  \ref{mainrem}, Theorems \ref{cl} and \ref{kirprop}. $\Box$

\end{document}